\magnification=\magstep1
\baselineskip=18pt
\font\bbd = msbm10 at 10pt
\font\bbds = msbm8 at 8pt
\font\bbds = msbm7 at 8pt
\font\csc = cmcsc10 at 10pt
\font\small = cmr8 at 8pt
\font\frak = eufm10 at 10pt
\font\fraksmall = eufm8 at 8pt
\font\xmplbx = cmbx10 scaled \magstep1
\outer\def\beginsection#1{\vskip 0pt plus.3\vsize\penalty-10
                          \vskip 0pt plus-.3\vsize\bigskip\medskip\vskip\parskip
                          \message{#1}\leftline{\bf #1}
                          \nobreak\smallskip}
\outer\def\proclaim #1 #2\par{\medbreak
      \noindent{\bf#1\enspace}{\sl#2}\par
      \ifdim\lastskip<\medskipamount \removelastskip\penalty55\medskip\fi}
\def\proof {\noindent {\it Proof.} \ \ }

\def\sqr#1#2{{\vcenter{\hrule height.#2pt
              \hbox{\vrule width.#2pt height#1pt \kern#1pt \vrule width.#2pt}
                    \hrule height.#2pt}}}
\def\qe{\mathchoice\sqr54\sqr54\sqr{2.1}3\sqr{1.5}3}

\def\qed{\quad\qe}

\def\C{\mathord{\hbox{\bbd C}}}
\def\c{\mathord{\hbox{\bbds C}}}

\def\F{\mathord{\hbox{\bbd F}}}

\def\V{\mathord{\hbox{\frak V}}}
\def\v{\mathord{\hbox{\fraksmall V}}}

\def\A{\mathord{\hbox{\frak a}}}

\def\F{\mathord{\hbox{\frak F}}}
\def\g{\mathord{\hbox{\fraksmall g}}}
\def\G{\mathord{\hbox{\frak g}}}
\def\n{\mathord{\hbox{\frak n}}}
\def\s{\mathord{\hbox{\frak s}}}
\def\NS{\mathord{\hbox{\frak N}}}
\def\ns{\mathord{\hbox{\fraksmall N}}}
\def\l{\mathord{\hbox{\fraksmall L}}}
\def\L{\mathord{\hbox{\frak L}}}

\def\Z{\mathord{\hbox{\bbd Z}}}
\def\z{\mathord{\hbox{\bbds Z}}}
\def\lsemi{\mathbin{\hbox{\bbd n}}}

\def\setbrackets #1 #2{\bigr\{\,#1\bigm|#2\,\bigr\}}

\def\dd{\mathord{\hbox{$\partial$}}}

\centerline { \xmplbx Conformal Modules}
\medskip
\centerline { by }
\medskip
\centerline {Shun-Jen Cheng\footnote{$^\dagger$}{\small partially supported  
by NSC grant
86-2115-M-006-012 of the ROC}
 and Victor G. Kac\footnote{$^{\dagger\dagger}$}{\small partially supported  
by NSF grant DMS-9622870}}
\medskip
\vbox{\small\baselineskip=10pt  $^\dagger$Department of Mathematics,
National Cheng-Kung University, Tainan,
Taiwan\hfill\break\indent\phantom{$^\dagger$}chengsj@mail.ncku.edu.tw}
\medskip
\vbox{\small\baselineskip=10pt $^{\dagger\dagger}$Department of Mathematics,  
MIT, Cambridge, MA 02139,
USA\hfill\break\indent\phantom{$^\dagger\dagger$}kac@math.mit.edu}
\vskip 0.5in
{\noindent\bf Abstract.}  In this paper we study a class of modules over
infinite-dimensional Lie (super)algebras, which we call conformal modules.   
In particular we classify and construct explicitly all irreducible conformal  
modules over the Virasoro and the $N{=}1$ Neveu-Schwarz algebras, and over the  
current algebras.

\beginsection{0. Introduction}

Conformal module is a basic tool for the construction of free field
realization of infinite-dimensional Lie (super)algebras in conformal field
theory.  This is one of the reasons to classify and construct such modules.   
In the present paper we solve this problem under the irreducibility
assumption for the Virasoro and the Neveu-Schwarz algebras, for the
current algebras and their semidirect sums.  Since complete reducibility
does not hold for conformal modules, one has to discuss the extension
problem.  This problem is solved in [1].

The basic idea of our approach is to use three (more or less) equivalent
languages.  The first is the language of local formal distributions, the
second is the language of modules over conformal algebras, and the third is  
the language of conformal modules over the annihilation subalgebras.  The
problem is solved using the third language by means of the crucial Lemma 3.1.
Note that conformal modules over Lie algebras of Cartan type were studied
in [6], where, in particular, a proof of Corollary 3.3 is contained.

This paper is organized as follows.  In Section 1 the concepts of a Lie
(super)algebra of formal distributions and of a conformal (super)algebra are
recalled.  They are more or less equivalent notions.  In Section 2 we
introduce conformal modules over a Lie (super)algebra of formal distributions  
and clarify their connections with modules over the corresponding conformal  
(super)algebra.  We then show that modules over a
conformal (super)algebra are in 1-1 correspondence with modules over the
annihilation subalgebra of the associated Lie (super)algebra of formal
distributions.  At the end of Section 2 examples of
conformal modules over the Virasoro, the current and Neveu-Schwarz algebras  
and their various semidirect sums are constructed.  In Section 3 we first
prove the key lemma (Lemma 3.1) and with its help classify all irreducible
finite conformal modules over the annihilation subalgebra of the
above-mentioned Lie (super)algebras.  The main result, stated in Theorem 3.2, 
which describes all finite irreducible modules over the conformal
(super)algebras in question (hence all irreducible finite conformal modules
over the corresponding Lie (super)algebras), is then immediate.

Roughly speaking, the main claim of the present paper is that all non-trivial
modules over current, Virasoro, $N{=}1$ superconformal algebras
and their semidirect sums are non-degenerate. For $N{>}1$ superconformal
algebras interesting degeneracies occur in
some non-trivial conformal modules. These effects are studied in [5].

\beginsection{1. Preliminaries on local formal distributions and conformal
superalgebras}

A {\it formal distribution} (usually called a field by physicists) with
coefficients in a complex vector space $U$ is a generating series of the
form:
$$a(z)=\sum_{n\in\z}a_{(n)}z^{-n-1},$$
where $a_{(n)}\in U$ and $z$ is an indeterminate.

Two formal distributions $a(z)$ and $b(z)$ with coefficients in a Lie
superalgebra $\G$ are called (mutually) {\it local} if for some $N\in\Z_+$
one has:
$$(z-w)^N[a(z),b(w)]=0.\eqno{(1.1)}$$
Introducing the {\it formal delta function}
$$\delta(z-w)=z^{-1}\sum_{n\in\z}({z\over w})^n,$$
we may write a condition equivalent to (1.1):
$$[a(z),b(w)]=\sum_{j=0}^N(a_{(j)}b)(w){{\dd_{w}^{(j)}}\delta(z-w)},\eqno{(1.2)}$$
(here $\dd^{(j)}_w$ stands for ${1\over {j!}}{{\partial^j}\over {\partial
w^j}}$) for some formal distributions $(a_{(j)}b)(w)$ ([3], Theorem 2.3),
which are uniquely determined by the formula
$$(a_{(j)}b)(w)={\rm Res}_z(z-w)^j[a(z),b(w)].\eqno{(1.3)}$$
Formula (1.3) defines a $\C$-bilinear product $a_{(j)}b$ for each $j\in\Z_+$  
on the space of all formal distributions with coefficients in $\G$.

Note also that the space (over $\C$) of all formal distributions with
coefficients in $\G$ is a (left) module over $\C[\dd_z]$, where the action of  
$\dd_z$ is defined in the obvious way, so that $\dd_z a(z)=\sum_n(\dd
a)_{(n)}z^{-n-1}$, where $(\dd a)_{(n)}=-na_{(n-1)}$.

The Lie superalgebra $\G$ is called a {\it  Lie superalgebra of formal
distributions} if there exists a family $\F$ of pairwise local formal
distributions whose coefficients span $\G$. In such a case we say that the
family $\F$ {\it spans} $\G$. We will write $(\G,\F)$ to emphasize the
dependence on $\F$.

The simplest example of a Lie superalgebra of formal distributions is the
{\it current superalgebra} $\tilde{\G}$ associated to a finite-dimensional
Lie superalgebra $\G$:
$$\tilde{\G}=\G\otimes\C[t,t^{-1}].$$
It is spanned by the following family of pairwise local formal distributions  
$a\in\G$:
$$a(z)=\sum_{n\in\z}(a\otimes t^n)z^{-n-1}.$$
Indeed, it is immediate to check that
$$[a(z),b(w)]=[a,b](w)\delta(z-w).$$

The simplest example beyond current algebras is the (centerless) {\it
Virasoro algebra}, the Lie algebra $\V$ with basis $L_n$ ($n\in\Z$) and
commutation relations
$$[L_m,L_n]=(m-n)L_{m+n}.$$
It is spanned by the local formal distribution $L(z)=\sum_{n\in\z}L_n
z^{-n-2}$, since one has:
$$[L(z),L(w)]=\dd_{w}L(w)\delta(z-w)+2L(w)\dd_{w}\delta(z-w).\eqno{(1.4)}$$
The next important example is the semidirect sum $\V\lsemi \tilde{\G}$ with  
the usual commutation relations between $\V$ and $\tilde{\G}$:
$$[L_m,a\otimes t^n]=-na\otimes t^{m+n},$$
which is equivalent to
$$[L(z),a(w)]=\dd_{w}a(w)\delta(z-w)+a(w)\dd_w\delta(z-w).$$

Given a Lie superalgebra of formal distributions $(\G,\F)$, we may always
include $\F$ in the minimal family $\F^c$ of pairwise local distributions
which is closed under $\dd$ and all products (1.3) ([3], Section 2.7).  Then  
$\F^c$
is a {\it conformal superalgebra} with respect to the products (1.3). Its
definition is given below [3]:

A {\it conformal superalgebra} is a left ($\Z_2$-graded) $\C[\dd]$-module
$R$ with a $\C$-bilinear product $a_{(n)}b$ for each $n\in\Z_+$ such that the  
following axioms hold ($a,b,c\in R; m,n\in\Z_+$ and $\dd^{(j)}={1\over
{j!}}\dd^j$):
\item{(C0)} \quad $a_{(n)}b=0$, for $n>>0$,
\item{(C1)} \quad $(\dd a)_{(n)}b=-na_{(n-1)}b,$
\item{(C2)} \quad $a_{(n)}b=(-1)^{{{\rm deg}a}{{\rm
deg}b}}\sum_{j=0}^{\infty}(-1)^{j+n+1}{\dd^{(j)}}(b_{(n+j)}a)$,
\item{(C3)} \quad $a_{(m)}(b_{(n)}c)=\sum_{j=0}^{\infty}{m\choose
j}(a_{(j)}b)_{(m+n-j)}c+(-1)^{{{\rm deg}a}{{\rm deg}b}}b_{(n)}(a_{(m)}c)$.

\noindent Of course, conformal algebra coincides with its even part, i.e.~${\rm
deg}a=0$ for all $a\in R$ in this case. Note also the following consequence of
(C1) and (C2):
\item{(C2')} \quad $a_{(n)}\dd b=\dd(a_{(n)}b)+na_{(n-1)}b$,

\noindent hence $\dd$ is a derivation of all products (1.3).

Conversely, assuming for simplicity (cf. Lemma 2.2b) that $R=\oplus_{i\in  
I}\C[\dd]a^i$ is a free (as a
$\C[\dd]$-module) conformal superalgebra, we may associate to $R$ a Lie
superalgebra of formal distributions $(\G(R),\F)$ with basis $a^i_{(m)}$
($i\in I$, $m\in\Z$) and $\F=\{a^i(z)=\sum_na^i_{(n)}z^{-n-1}\}_{i\in I}$
with bracket (cf.~(1.2)):
$$[a^i(z),a^j(w)]=\sum_{k\in\z_+}(a^i_{(k)}a^j)(w){{\dd_w^{(k)}}}\delta(z-w),\eqno{(1.5)}$$
so that $\F^c=R$.

Formula (1.5) is equivalent to the following commutation relations
($m,n\in\Z$; $i,j\in I$):
$$[a^i_{(m)},a^j_{(n)}]=\sum_{k\in\z_+}{m\choose
k}(a^i_{(k)}a^j)_{(m+n-k)}.\eqno{(1.6)}$$
It follows that the $\C$-span of all $a^i_{(n)}$ with $n\in\Z_+$ is a
subalgebra of the Lie superalgebra $\G(R)$.  We denote this subalgebra by
$\G(R)_+$ and call it the {\it annihilation subalgebra}.

For example $\V_+=\sum_{j\ge -1}\C L_j$ and $\tilde{\G}_+=\G\otimes \C[t]$.

The simplest examples of a conformal superalgebra is the {\it current conformal  
superalgebra} associated to a finite-dimensional Lie superalgebra $\G$:
$$R(\tilde{\G})=\C[\dd]\otimes_{\c}\G,$$
with the products defined by:
$$a_{(0)}b=[a,b],\quad a_{(j)}b=0,\quad {\rm for\ }j>0,\ a,b\in\G,$$
and the {\it Virasoro conformal algebra} $R(\V)=\C[\dd]\otimes_{\c} L$ with  
products (cf.~(1.4)):
$$L_{(0)}L=\dd L,\quad L_{(1)}L=2L,\quad L_{(j)}L=0,\quad {\rm for\ }j>1.$$

Their semidirect sum is $R(\V)\lsemi R(\tilde{\G})$ with additional non-zero
products $L_{(0)}a=\dd a$ and $L_{(1)}a=a$, for $a\in\G$. These examples
are the conformal superalgebras associated to the Lie superalgebras of
formal distributions described above.

The simplest superextension of the Virasoro algebra is the well-known
(centerless) Neveu-Schwarz algebra $\NS$ which, apart from even basis
elements $L_n$, has odd basis elements $G_r$, $r\in{1\over 2}+\Z$, with
commutation  relations:
$$[G_r,L_n]=(r-{n\over 2})G_{r+n},\quad [G_r,G_s]=2L_{r+s}.$$
The corresponding annihilation subalgebra in this case is $\NS_+=\sum_{n\ge  
-1}\C L_n+\sum_{r\ge -{1\over 2}}\C G_r$.  The conformal superalgebra,
associated to $\NS$, is $R(\NS)=\C[\dd]\otimes_{\c}L+\C[\dd]\otimes_{\c}G$
with additional non-zero products:
$$L_{(0)}G=\dd G,\quad G_{(0)}L={1\over 2}\dd G, \quad
L_{(1)}G=G_{(1)}L={3\over 2}G,\quad G_{(0)}G=2L.$$

Other examples treated in this paper are {\it supercurrent algebras}
$$\tilde{\G}_{\rm super}=\G\otimes_{\c}\C[t,t^{-1},\theta],$$
where $\theta$ is an odd indeterminate.  The generating family $\F$ of
pairwise local formal distributions consists of currents $a(z)$ ($a\in\G$)
introduced above and supercurrents
$${a}^{\theta}(z)=\sum_{n\in\z}(a\otimes t^n\theta)z^{-n-1}.$$
Of course its associated conformal superalgebra is $R(\tilde{\G}_{\rm
super})=\C[\dd]\otimes(\G\oplus \G^{\theta})$, where $\G^{\theta}$ is an
identical copy of $\G$, but with reversed parity.  $R(\tilde{\G}_{\rm
super})$ extends $R(\tilde{\G})$ by the additional non-trivial product
$$a_{(0)}b^{\theta}=[a,b]^{\theta}, \quad a,b\in\G.$$

The final example treated in this paper is the semidirect sum
$\NS\lsemi\tilde{\G}_{\rm super}$, which is defined by letting
$L_n=-t^n(t{\partial\over{\partial t}}+{{n+1}\over 2}\theta{\partial\over
{\partial\theta}})$ and $G_r=-t^{r+{1\over 2}}(\theta{\partial\over{\partial  
t}}-{\partial\over{\partial\theta}})$ for $n\in\Z$ and $r\in{1\over 2}+\Z$.   
Its corresponding conformal superalgebra $R(\NS\lsemi\tilde{\G}_{\rm
super})= R(\NS)\lsemi R(\tilde{\G}_{\rm super})$ has the following
additional non-trivial
products:
$$L_{(0)}a^{\theta}=\dd a^{\theta},\quad L_{(1)}a^{\theta}={1\over
2}a^{\theta},\quad G_{(0)}a^{\theta}=a,\quad G_{(0)}a=\dd a^{\theta}, \quad  
G_{(1)}a= a^{\theta}.$$

\beginsection{2. Preliminaries on conformal modules}

Let $(\G,\F)$ be a Lie superalgebra of formal distributions, and let $V$ be  
a $\G$-module.  We say that a formal distribution $a(z)\in\F$ and a formal
distribution $v(z)=\sum_{n\in\z}v_{(n)}z^{-n-1}$ with coefficients in $V$ are  
{\it local} if
$$(z-w)^N a(z)v(w)=0,\quad {\rm for\ some\ }N\in\Z_{+}.\eqno{(2.1)}$$
It follows from [3] {Section  2.3} that (2.1) is equivalent to
$$a(z)v(w)=\sum_{j=0}^{N-1}(a_{(j)}v)(w){\partial_{w}^{(j)}}\delta(z-w),\eqno{(2.2)}$$
for some formal distributions $(a_{(j)}v)(w)$ with coefficients in $V$,
which
are uniquely determined by the formula
$$(a_{(j)}v)(w)={\rm Res}_{z}(z-w)^ja(z)v(w).$$
Formula (2.2) is obviously equivalent to
$$a_{(m)}v_{(n)}=\sum_{j\in\z_+}{m\choose j}(a_{(j)}v)_{(m+n-j)}.\eqno{(2.3)}$$

\noindent {\it Example 2.1.} Consider the following representation of the
(centerless) Virasoro algebra in the vector space $V$ with basis $v_{(n)}$,  
$n\in\Z$, over $\C$
$$L_mv_{(n)}=((\Delta-1)(m+1)-n)v_{(m+n)}+\alpha v_{(m+n+1)},$$
where $\Delta,\alpha\in\C$.  In terms of formal distributions $L(z)$ and
$v(z)$ this can be written as follows:
$$L(z)v(w)=(\dd+\alpha)v(w)\delta(z-w)+\Delta v(w)\dd_{w}\delta(z-w).\eqno{(2.4)}$$
Hence $L(z)$ and $v(z)$ are local.

Suppose that $V$ is spanned over $\C$ by the coefficients of a family $E$ of  
formal distributions such that all $a(z)\in\F$ are local with respect to all  
$v(z)\in E$.  Then we call $(V,E)$ a {\it conformal module over $(\G,\F)$}.

The following is a representation-theoretic analogue (and a generalization)  
of Dong's lemma (see [3], Section  3.2).

\proclaim{Lemma 2.1.} {\rm [4]} Let $V$ be a module over a Lie superalgebra  
$\G$, let $a(z)$ and $b(z)$ (respectively $v(z)$) be formal distributions
with coefficients in $\G$ (respectively in $V$).  Suppose that all pairs
$(a,b)$, $(a,v)$ and $(b,v)$ are local.  Then the pairs $(a_{(j)}b,v)$
and $(a,b_{(j)}v)$ are local for all $j\in\Z_+$.

This lemma shows that the family $E$ of a conformal module $(V,E)$ over
$(\G,\F)$ can always be included in a larger family $E^c$ which is still
local with respect to $\F$, hence to $\F^c$, and such that $\dd E^c\subset E^c$ and
$a_{(j)}E^c\subset E^c$ for all $a\in\F$ and $j\in\Z_+$.

It is straightforward to check the following properties for $a,b\in\F$ and
$v\in E^c$:
$$[a_{(m)},b_{(n)}]v=\sum_{j=0}^m {m\choose j}(a_{(j)}b)_{(m+n-j)}v,\eqno{(2.5)}$$
$$(\dd a)_{(n)} v=[\dd,a_{(n)}]v=-na_{(n-1)}v.\eqno{(2.6)}$$
(Here $[\cdot,\cdot]$ is the bracket of operators on $E^c$.) It follows from  
(2.5) (by induction on m) and (2.6) that $a_{(j)}E^c\subset E^c$ for all
$a\in\F^c$ and $j\in\Z_+$.

Thus any conformal module $(V,E)$ over a Lie superalgebra of formal
distributions $(\G,\F)$ gives rise to a module $M=E^c$ over the conformal
superalgebra $R=\F^c$, defined as follows.  It is a (left) $\Z_2$-graded
$\C[\dd]$-module equipped with a family of $\C$-linear maps $a\rightarrow
a_{(n)}^M$ of $R$ to ${\rm End}_{\c}M$, for each $n\in\Z_+$, such that the
following properties hold (cf.~(2.5) and (2.6)) for $a,b\in R$ and
$m,n\in\Z_+$ :
\item{(M0)} $\quad a^M_{(n)}v=0$, for $v\in M$ and $n>>0$,
\item{(M1)} $\quad [a^M_{(m)},b^M_{(n)}]=\sum_{j=0}^m{m\choose
j}(a_{(j)}b)^M_{(m+n-j)}$,
\item{(M2)} $\quad (\dd a)^M_{(n)}=[\dd,a^M_{(n)}]=-na_{(n-1)}^M$.

Conversely, suppose that a conformal superalgebra $R=\oplus_{i\in
I}\C[\dd]a^i$ is a free
$\C[\dd]$-module and consider the associated Lie superalgebra of formal
distributions $(\G(R),\F)$ (see Section  1).  Let $M$ be a module over the
conformal superalgebra $R$ and suppose (cf. Lemma 2.2b) that $M$ is a free  
$\C[\dd]$-module
with $\C[\dd]$-basis $\{v^{\alpha}\}_{\alpha\in J}$. This gives rise to
a conformal  $\G(R)$-module $V(M)$ with basis $v_{(n)}^{\alpha}$, where $\alpha\in J$ and $n\in\Z$, defined by (cf.~(2.2)):
$$a^{i}(z)v^{\alpha}(w)=\sum_{j\in\z_+}(a^{i}_{(j)}v^{\alpha})(w){{\dd_{w}^{(j)}}}\delta(z-w).\eqno{(2.7)}$$

A conformal module $(V,E)$ (respectively module $M$) over a Lie superalgebra  
of formal distributions $(\G,\F)$ (respectively over a conformal superalgebra  
$R$) is called {\it finite}, if $E^c$ (respectively $M$) is a finitely
generated $\C[\dd]$-module.

A conformal module $(V,E)$ over $(\G,\F)$ is called {\it irreducible} if
there is no non-trivial invariant subspace which contains all $v_{(n)}$,
$n\in\Z$, for some non-zero $v\in E^c$.  Clearly a conformal module is  
irreducible if
and only if the associated module $E^c$ over the conformal superalgebra
$\F^c$ is irreducible (in the obvious sense).

The above discussions, combined with the following lemma, reduce the
classification of finite conformal modules over a Lie superalgebra of formal  
distributions $(\G,\F)$ to the classification of finite modules over the
corresponding conformal superalgebra.

\proclaim{Lemma 2.2.} {\rm [4]}
(a) Let $M$ be a module
over a conformal superalgebra $R$ and
let $v\in M$ be such that $\dd v=\lambda v$ for some $\lambda\in\C$.  Then
$v$ is an invariant, i.e.~$R_{(m)}v=0$ for all $m\in\Z_+$.
\item{(b)} Suppose that $M$ is a finite module over a conformal
superalgebra and suppose that $M$ has no non-zero invariants. Then $M$ is a  
free $\C[\dd]$-module.

\noindent{\it Remark 2.1.} Given a module $M$ over a conformal superalgebra  
$R$, we may change its structure as a $\C[\dd]$-module by replacing $\dd$ by  
$\dd+A$, where $A$ is an endomorphism over $\C$ of $M$ which commutes with
all $a_{(n)}^M$ (this will not affect axiom (M2)).

Note that the maps $a\rightarrow a_{(n)}$ of $R$ to ${\rm End}_{\c}R$
defines an $R$-module, called the {\it adjoint module}.

Example 2.1 gives a $2$-parameter family of (irreducible) modules over the
Virasoro conformal algebra.  Note also that the well-known family of graded  
Virasoro modules given by $L_mv_{(n)}=((\Delta-1)(m+1)-n+\alpha)v_{(m+n)}$ is  
conformal, but is finite if and only if $\alpha=0$.

The following simple observation, which follows from definitions,
is fundamental for representation theory of
conformal superalgebras.

\proclaim{Proposition 2.1.} Consider the Lie superalgebra of formal
distributions $(\G(R),\F)$ defined by (1.5) and let $\G(R)_+$ be the
annihilation subalgebra of $\G(R)$.  Denote by $\G(R)^+$ the semidirect
product of the $1$-dimensional Lie subalgebra $\C\dd$ and the ideal $\G(R)_+$  
with the action of $\dd$ on $\G(R)_+$ given by
$\dd(a^i_{(n)})=-na^i_{(n-1)}$.  Then a module $M$ over the conformal
superalgebra $R$ is precisely a $\G(R)^+$-module $M$ (over $\C$) such that
$$a^i_{(n)}v=0, \quad {\rm for\ each\ } v\in M {\rm \ and\ }n>>0.\eqno{(2.8)}$$

\proclaim{Corollary 2.1.} Let $R=\oplus_{i\in I}\C[\dd]a^i$ be a conformal
superalgebra and $M=\oplus_{j\in J}\C[\dd]v^j$ be a free $\C[\dd]$-module.   
Then, given $a^i_{(n)}v^j\in M$ for all $i\in I$, $j\in J$, $n\in\Z_+$, which  
is $0$ for $n>>0$, we may extend uniquely the action of $a^i_{(n)}$ to the
all of $R$ on $M$ using (M2).  Suppose that (M1) holds for all $a=a^i_{(m)}$,  
$b=a^j_{(n)}$.  Then $M$ is an $R$-module.

Using Proposition 2.1 and Corollary 2.1, one can construct large families of  
finite modules over conformal superalgebras, hence corresponding modules
over Lie superalgebras of formal distributions.

In conclusion we will list more examples of modules over conformal
superalgebras.  In Section 3
the irreducible ones listed below will turn out to exhaust the list of all
irreducible finite conformal modules over these conformal superalgebras.

\noindent {\it Example 2.2.} Let $\V_{0}=\sum_{j\ge 0}\C L_j$ and consider a  
representation $\pi$ of the Lie algebra $\V_{0}$ in a finite-dimensional
(over $\C$)
vector space $U$.  Let $A$ be an endomorphism of $U$ commuting with all
$\pi(L_j)$ ($j\in\Z_+$).  Then $\C[\dd]\otimes U$ is a finite module over the  
conformal algebra $\V$ defined by the following formulas ($u\in U$):
$$L_{(0)}u=(\dd +A)u,\quad L_{(j)}u=\pi(L_{j-1})u,{\ \rm for\ }j\ge 1.$$
For example we can take $\pi(L_0)=B$, where $B$ is an endomorphism of $U$
commuting with $A$.  Then
$$L_{(0)}u=\dd u+Au,\quad L_{(1)}u=Bu,\quad L_{(j)}u=0,\ j\ge 1,$$
defines a finite module over the centerless conformal algebra $\V$, which we  
denote by $M_{\v}(A,B)$.  Taking ${\rm dim} U=1$, $A=\alpha$ and $B=\Delta$,  
where $\alpha,\Delta\in\C$, gives Example 2.1.

Translating back to the language of Lie algebras of formal distributions,
Example 2.2 gives the following family of finite conformal modules over the  
Virasoro algebra in the space $U\otimes\C[t,t^{-1}]$ (we let, as usual,
$u_{(n)}=u\otimes t^{n}$):
$$L_m
u_{(n)}=(Au)_{(m+n+1)}-(m+n+1)u_{(m+n)}+\sum_{j=0}^{\infty}{{m+1}\choose
{j+1}}(\pi(L_j)u)_{(m+n-j)}.$$
The above special case defined by a pair $(A,B)$ of the commuting operators  
in $U$ is given by:
$$L_{m}u_{(n)}=(Au)_{(m+n+1)}+((m+1)Bu-(m+n+1)u)_{(m+n)}.$$
We keep the notation $M_{\v}(A,B)$ for this $\V$-module.  Note that the
module $M_{\v}(A,B)$ is irreducible if and only if the $\V_{0}$-module is
irreducible and non-trivial, i.e.~if and only if ${\rm dim}U=1$ and
$B=\Delta\not=0$.  We denote these $\V$-modules again by
$M_{\v}(\alpha,\Delta)$, $\alpha,\Delta\in\C$ and $\Delta\not=0$.

\noindent{\it Remark 2.2.} (cf.~Remark 2.1 and Example 2.2.) Suppose that
the annihilation subalgebra $\G_+$ of the Lie algebra
$\G^+=\C[\dd]\lsemi\G_+$, contains an element $L_{-1}$ such that ${\rm
ad}L_{-1}=\dd$ on $\G_+$.  Then $\G^+$ is a direct sum of ideals
$\C(\dd-L_{-1})$ and $\G_+$. Hence in this case the study of conformal
modules reduces to the study of modules over $\G_+$ satisfying (2.8).  This  
remark applies to all cases except for the current and the supercurrent
algebras.

\noindent{\it Example 2.3.} Consider the current Lie superalgebra
$\tilde{\G}$ and the associated conformal superalgebra $R(\tilde{\G})$. Let  
$\pi$ be a representation of $\tilde{\G}_+=\G\otimes_{\c}\C[t]$ in a
finite-dimensional vector space $U$, such that $(t^n\otimes\G)U=0$ for
$n>>0$.  This defines on the space $U\otimes\C[t,t^{-1}]$ the structure of a  
conformal module over $\tilde{\G}$ by the formula:
$$(a\otimes t^m)u_{(n)}=\sum_{j\in \z_+}{m\choose j}(\pi(a\otimes
t^j)u)_{(m+n-j)}.$$

A special case of this construction is to take a finite-dimensional
representation $\pi$ of the Lie superalgebra $\G$ in a finite-dimensional
vector space $U$ and extend it to $\tilde{\G}_+$ by letting  $\G\otimes
t\C[t]$ act trivially.  Then we have $(a\otimes
t^m)u_{(n)}=(\pi(a)u)_{(m+n)}$.  Translating back to the language of modules  
over the conformal algebra $R(\tilde{\G})$ we obtain the free  
$\C[\dd]$-module $\C[\dd]\otimes_{\c}U$ with
$$a_{(0)}u=\pi(a)u,\quad a_{(n)}u=0\quad{\rm for\ } n>0,{\ \rm where\ }u\in U$$
We will denote this module by $M_{\tilde{\g}}(\pi)$.  It is irreducible if  
and only if $\pi$ is irreducible.  In this case we will denote the module by  
$M_{\tilde{\g}}(\Lambda)$, where $\Lambda$ is the highest weight of $U$.

\noindent{\it Example 2.4.} Consider the $N{=}1$ Neveu-Schwarz algebra $\NS$  
with associated conformal superalgebra $R(\NS)$.  Let $\NS_+$ denote the
corresponding annihilation subalgebra.  Consider a finite-dimensional
representation $(\pi,U)$ of $\NS_{0}$, the subalgebra of
$\NS_+$ spanned by elements of non-negative modes.  Let
$U^{\theta}$ denote an identical copy
of $U$ with reversed parity.  Then the following gives a structure of a
module over $R(\NS)$ on the free $\C[\dd]$-module $\C[\dd]\otimes
(U\oplus U^{\theta})$:
$$L_{(0)}u=(\dd +A)u,\quad L_{(j)}u=\pi(L_{j-1})u,\quad
L_{(0)}u^{\theta}=(\dd+A)u^{\theta},$$  
$$L_{(1)}u^{\theta}=(\pi(L_{0})u)^{\theta}+{1\over 2}u^{\theta},\quad
L_{(j+1)}u^{\theta}=(\pi(L_{j})u)^{\theta}+{{j+1}\over 2}\pi(G_{j-{1\over  
2}})u,\quad
G_{(0)}u=u^{\theta},$$
$$G_{(j)}u=\pi(G_{j-{1\over 2}})u,\quad
G_{(0)}u^{\theta}=(\dd u+Au),\quad
G_{(j)}u^{\theta}=(\pi(G_{j-{1\over 2}})u)^{\theta}+2\pi(L_{j-1})u,$$
where $j\ge 1$, $u\in U$, and $A$ is an operator acting on $U$, commuting with all
$\pi(\NS_0)$.  In particular let $U=\C u$ be the one-dimensional  
$\NS_0$-module with action $L_{0}u=\Delta u$, with
$0\not=\Delta\in\C$, all other generators acting trivially.  Then for  
$A=\alpha\in\C$ arbitrary, we obtain an irreducible
module over $R(\NS)$.  Hence as in the case of the Virasoro algebra we get a  
$2$-parameter family of finite irreducible conformal modules.  Denote this
family by $M_{\ns}(\alpha,\Delta)$.

\noindent{\it Example 2.5.} Consider the supercurrent algebra $\tilde{\G}_{\rm
super}$.  Let $R(\tilde{\G}_{\rm super})$, ${\tilde{\G}_{\rm super+}}$ be as  
usual.  Let $(\pi,U)$ be a finite-dimensional representation of
${\tilde{\G}_{\rm super+}}$.  We obtain a module over $R(\tilde{\G}_{\rm
super})$ as in the case of current algebra by setting for $n\ge 0$:
$$a_{(n)}u=\pi(a\otimes t^n)u,\quad {a}^{\theta}_{(n)}u=\pi(a\otimes
t^n\theta)u,\quad a\in\G,\ u\in U.$$  Denote these modules by
$M_{\tilde{\g}_{\rm super}}(\pi)$.
In the special case when $U$ is a finite-dimensional irreducible
representation of $\G$ of highest weight $\Lambda\not=0$, extended to
${\tilde{\G}_{\rm super+}}$ trivially, the associated module over
$R(\tilde{\G}_{\rm super})$ is irreducible and finite.  We denote this
module by $M_{\tilde{\g}_{\rm super}}(\Lambda)$.

\beginsection{3. The key lemma and classification of finite irreducible
conformal modules}

Let $(\G,\F)$ be a Lie superalgebra of formal distributions.  For each
$N\in\Z_+$ let
$$\G_N=\sum_{a\in\F,n\ge N}\C a_{(n)}.$$
Suppose that $(\G,\F)$ is {\it regular} [3], i.e.~there exists a derivation  
$\dd$ of $\G$ such that $\dd(a_{(n)})=-na_{(n-1)}$ for all $a\in\F$ and
$n\in\Z$.  Obviously, $(\G(R),\F)$, where $R$ is a conformal superalgebra,
is regular, hence all examples considered above
are regular.  Then $\G_+=\G_0$ is the annihilation subalgebra, which is
$\dd$-invariant and, due to Proposition 2.1, we have to study representations  
of the Lie superalgebra $\G^+=\C\dd \lsemi \G_+$, called the {\it extended  
annihilation subalgebra}.  This leads us to consider the
following more abstract situation.

Let $\L$ be a Lie superalgebra over $\C$ with a distinguished element $\dd$  
and a descending sequence of
subspaces
$\L\supset\L_{0}\supset\L_{1}\supset\L_{2}\supset\cdots\supset\L_{n}
\supset\cdots$, such that $[\dd,\L_{k}]=\L_{k-1}$, for all $k>0$. Let $V$ be an
$\L$-module such that given any
$v\in V$ there exists a non-negative integer $N$ (depending on $v$)  such that
$\L_{N}v=0$.  We will call such
$\L$-modules {\it conformal $\L$-modules}.  A conformal $\L$-module is
called {\it finite} if it is finitely generated as a $\C[\dd]$-module.  Our  
goal is to describe irreducible finite conformal $\L$-modules.

Let $V$ be a conformal $\L$-module. Let $V_n=\{v\in V|\L_nv=0\}$, and let
$N$ be the minimal non-negative integer such that $V_N\not=0$ (it
exists by definition).  Let $U=V_N$.  Now we can state
the key lemma.

\proclaim{Lemma 3.1.} Suppose that $N\ge 1$.  Then  
$\C[\dd]U=\C[\dd]\otimes_{\c}U$ and hence $\C[\dd]U\cap U=U$.
In particular $U$ is a finite-dimensional vector
space if $V$ is finite.

\proof  Let ${\rm L}_a$ and ${\rm R}_a$ denote the left and right multiplication
by the element $a$, respectively.  Using ${\rm R}_a={\rm L}_a-{\rm ad}a$ and  
the binomial formula, we
get the following well-known formula in any associative algebra $A$,
$$ga^{k}=\sum_{j=0}^k{k\choose j}a^{k-j}(-{\rm ad}a)^jg,\quad a,g\in A.$$

Let $\{v_i\}$, $i\in I$, be a
$\C$-linearly independent set of vectors in
$U$ generating the $\C[\dd]$-module $\C[\dd]U$.
Suppose that $v=\sum_{i=1}^{n}p_i(\dd)v_i=0$, where
$p_i(\dd)\in\C[\dd]$, and not all $p_i(\dd)=0$.
Let $m$ be the maximal degree of the $p_i(\dd)$'s.  We write
$p_i(\dd)=\sum_{j=0}^mc_{ij}\dd^j$, $c_{ij}\in\C$, so that we have $c_{im}\not=0$,
for some $i$.  We have, since $N\ge 1$:
$$\eqalignno{\L_{N+m-1}\dd^k=&\sum_{j=0}^k{k\choose j}\dd^{k-j}({\rm
ad}\dd)^j(\L_{N+m-1})\cr =&\sum_{j=0}^k{k\choose j}\dd^{k-j}\L_{N+m-1-j}.\cr}$$
We have therefore
$$0=\L_{N+m-1}v=\sum_{i=1}^n c_{im}\L_{N-1}v_i=\L_{N-1}(\sum_{i=1}^nc_{im}v_i).$$
Since $\sum_{i=1}^nc_{im}v_i\not=0$, this contradicts the minimality of $N$.  
 Hence all the $p_i(\dd)=0$, proving the lemma.
$\qed$

\proclaim{Theorem 3.1.} Assume that $[\L_0,\L_i]\subset\L_i$ for all $i\ge 0$ and
that $\L=\C \dd+\C [\dd,\dd]+\L_0$. Let $V$ be a non-trivial irreducible
conformal $\L$-module.
Then either
$V\cong{\rm Ind}_{\L_{0}}^{\L}U$, where $U$ is an
irreducible $\L_{0}$-module (and ${\rm dim}U<\infty$ if $V$ is finite), or
else, provided that $\L_0$ is an ideal of $\L$, $V$ can be a (non-trivial)
finite-dimensional irreducible $\L/\L_0$-module.

\proof We continue to use the notation above.  Clearly $U$ is $\L_0$-invariant. 
First assume that $N\ge 1$.  By Lemma 3.1 we see that $V$ contains an
$\L$-submodule isomorphic to ${\rm Ind}_{\L_{0}}^{\L}U$.  Hence $V\cong {\rm
Ind}_{\L_{0}}^{\L}U$.  Clearly
$U$ must be irreducible over $\L_{0}$ in order for $V$ to be irreducible.
Conversely the same argument as in the proof of Lemma 3.1 shows that by applying
$\L_{k}$, for a suitable $k$, to a non-zero element of the form
$\sum_{i=1}^{n}q_i(\dd)v_i$, where $q_i(\dd)\in\C[\dd]$ and $v_i\in U$, one obtains
a non-zero element in $U$.  This implies that such induced modules are
irreducible.

Now suppose that $N=0$.  Let $u$ be a non-zero vector in $U$.  Then  
$U(\L)u=\C[\dd]U(\L_0)u=\C[\dd]u$, and hence $V=\C[\dd]u$ by irreducibility  
of $V$.  We consider two cases:

{\csc Case 1.} Suppose that $\L_0$ is an ideal of $\L$.  Then $\L_0$ acts
trivially on $V$.  Thus
$V$ is an irreducible $\C[\dd]$-module.  It follows that if
$\dd$ is even, then $V$ is 1-dimensional, and if $\dd$ is odd with
$[\dd,\dd]\not=0$, $V$ is either the trivial or a 2-dimensional $\L/\L_0$-module
with $[\dd,\dd]$ acting as a non-zero scalar.

{\csc Case 2.} Suppose that $\L_0$ is not an ideal of $\L$.  Then there  
exist $l_0,l_0'\in\L_0$ such that $[l_0,\dd]=\dd+l_0'$.  An easy induction  
argument shows that, $l_0\dd^iu=i\dd^iu+f_i(\dd)u$ for $i\ge 1$, where  
$f_i(\dd)$ is a polynomial in $\dd$ of degree at most $i-1$ with zero  
constant term.
Now if $v=\sum_{i=0}^na_i\dd^i u$ is a non-zero vector in $\C[\dd]\otimes \C  
u$ with $a_0\not=0$, then $v-({1\over n})l_0v$ is a non-zero vector of  
degree at most $n-1$.  Thus proceeding this way we see that $u$ is contained  
in the module generated by $v$. Therefore $\dd\C[\dd]\otimes\C u$ is the  
unique maximal submodule of $\C[\dd]\otimes \C u$ and hence  $V$ is the  
trivial module.  $\qed$

We will now apply the theorem above to classify irreducible conformal
modules over the Virasoro algebra, the current
algebra and their semidirect product.  In addition similar results can be
obtained for the corresponding $N{=}1$ extended superalgebras by slightly
modifying the arguments.  In order to do so, the following lemma is useful.

\proclaim{Lemma 3.2.} {\rm [5]} Let $\G$ be a Lie superalgebra and let $\n$ be
an ideal of $\G$. Assume that any finite-dimensional quotient of $\n$ is
solvable. Let $\A$ be an even subalgebra of
$\G$ such that $\n$ is a completely reducible ${\rm ad}\A$-module with no
trivial summand.  
Then $\n$ annihilates a non-zero vector in any finite-dimensional $\G$-module $V$.  In particular $\n$ acts trivially in any irreducible
finite-dimensional $\G$-module.

\proof First note that if $[a,b]=b$, then $b$ is nilpotent on any
finite-dimenisonal representation $V$.  To see this, note that if $v\in V$ is  
an eigenvector of $a$ with eigenvalue $\lambda$, then the condition
$[a,b]=b$ implies that $bv$ is an eigenvector of $a$ with eigenvalue
$\lambda+1$.  Thus there exists a non-zero $w\in V$ such that $bw=0$.  Let
$W$ be the space annihilated by $b$.  Then $W$ is $a$-invariant.  Since
$W\not=0$, it follows by induction on the dimension of $V$ that $b$ is
nilpotent on $V/W$.  Thus $b$ is nilpotent on $V$.

Taking the image of $\G$ in ${\rm End}_{\c}(V)$, we may assume that ${\rm
dim}_{\c}\G<\infty$, hence ${\rm dim}_{\c}{\A}<\infty$ and $\n$ is a
finite-dimensional solvable Lie superalgebra . Let $\s$ be an
irreducible ${\rm ad}\A$-submodule of $\n$. By the assumption, it is a
module with a non-zero highest weight.   
Hence there
exists $a\in\A$ and $b\in\s$ such that $[a,b]=b$, therefore $b$
is nilpotent on $V$.  Moreover all elements from the orbit ${\rm Ad}A\cdot
b$, where $A$ is the connected Lie group with Lie algebra $\A$, are
nilpotent on $V$.  Since $\s$ is ${\rm Ad}A$-irreducible, this orbit
spans $\s$, hence $\s$ is spanned by elements that are
nilpotent on $V$.  Thus any $\A$-submodule of  $\n$ is spanned by elements that are nilpotent on
$V$.

To prove the lemma, we may assume that $V$ is a faithful irreducible
$\G$-module.  Suppose
the lemma is not true, i.e.~$\n$ is non-zero.  Let $\n^{(i)}$ be the last
non-zero member of its derived series. By the above, $\n^{(i)}$ is spanned by
mutually commuting elements that are nilpotent on $V$, and hence $\n^{(i)}$  
annihilates a non-zero vector in $V$.
But $\n^{(i)}$ is an ideal of $\G$ and hence the subspace of $V$,
annihilated by $\n^{(i)}$, is a $\G$-submodule of $V$.  Thus $\n^{(i)}$
annihilates $V$ and so $V$ is not faithful, which is a contradiction. $\qed$

We are now in a position to classify finite conformal modules over the
Virasoro, Neveu-Schwarz, current and the supercurrent algebras and their
semidirect sums.  Due to Section 2 we only need to classify finite
modules over the corresponding (extended) annihilation subalgebras.  For  
each of these Lie
superalgebras of formal distributions, the corresponding annihilation
subalgebras are of course the corresponding subalgebras defined on the line,  
instead of the circle.  So we will use terminology like current algebras on  
the line etc.~to denote the corresponding annihilation subalgebras.

The following corollaries are immediate by Theorem 3.1 and Lemma 3.2.
\proclaim{Corollary 3.1.} Let $\G$ be a direct sum of finite-dimensional simple
Lie superalgebras (we allow commutative summands).  Let $\L=\C{d\over
{dt}}\lsemi \G[t]$, where
$\G[t]=\G\otimes\C[t]$ is the current algebra on the line.  Then every
non-trivial irreducible conformal module of
$\L$ is of the form ${\rm Ind}_{\g[t]}^{\l}U$, where $U$ is a finite
dimensional non-trivial irreducible $\G$-module or else it is the trivial
$\G[t]$-module on which ${d\over {dt}}$ acts as a non-zero scalar.

\proclaim{Corollary 3.2.} Let $\G$ be a direct sum of finite-dimensional simple
Lie superalgebras.  Let $\L=\C{d\over
{dt}}\lsemi \G[t,\theta]$, where $\G[t,\theta]=\G\otimes\C[t,\theta]$.  Then every
non-trivial irreducible conformal module of
$\L$ is of the form ${\rm Ind}_{\g[t,\theta]}^{\l}U$, where $U$ is a finite
dimensional non-trivial irreducible $\G$-module or else it is the trivial
$\G[t,\theta]$-module on which ${d\over {dt}}$ acts as a non-zero scalar.

\proclaim{Corollary 3.3.} Let $\V_+=\sum_{k\ge -1}\C{t^{k+1}{d\over{dt}}}$ be the
Virasoro algebra on the line.  Let
$L_i=-{t^{i+1}{d\over {dt}}}$ and let $\V_{0}=\sum_{k\ge 0}\C L_k$. Then
any non-trivial irreducible conformal
module of $\V_+$ is of the form ${\rm Ind}_{\v_{0}}^{\v_+}\C_{\lambda}$, where  
$\C_{\lambda}$ is a
non-trivial one-dimensional irreducible representation of $\V_{0}$, on which  
$L_0$ acts as $\lambda\in\C^*$ and $L_k$ act as $0$ for all
$k>0$.

\proclaim{Corollary 3.4.} Let $\L=\V_+\lsemi\G[t]$ such that $[L_k,a\otimes  
t^{n}]=-na\otimes t^{n+k}$, $a\in\G$.  Let $\L_0=\V_{0}\lsemi\G[t]$.  Then
every non-trivial irreducible $\L$-module
is of the form ${\rm Ind}_{\l_0}^{\l}U$, where $U$ is a non-trivial irreducible
$(\G\oplus\C L_{0})$-module with $\G[t]t$ and $L_k$, $k>0$, acting trivially.

\noindent{\it Remark 3.1.} Translating the modules over the annihilation
subalgebra $\V_+\lsemi\G[t]$ of Corollary 3.4 into the language of conformal  
modules we obtain a $3$-parameter family of non-trivial conformal modules
over $R(\V)\lsemi R(\tilde{\G})$.  We will denote these modules by
$M_{\v\lsemi\tilde{\g}}(\Lambda,\alpha,\Delta)$, where $\Lambda$ stands for  
the irreducible finite-dimensional $\G$-module of highest weight $\Lambda$.   
Clearly when $\Lambda\not=0$, $M_{\v\lsemi\tilde{\g}}(\Lambda,\alpha,\Delta)$  
is irreducible.  When $\Lambda=0$,
$M_{\v\lsemi\tilde{\g}}(0,\alpha,\Delta)$, for $\Delta\not=0$, is
irreducible.

For the Neveu-Schwarz algebra on the line we have the following description:

\proclaim{Corollary 3.5.}  Let $\NS_{+}=\sum_{n\ge -1}\C L_n
+\sum_{r\ge-{1\over 2}}\C G_r$ be the Neveu-Schwarz algebra on the line. Let  
$\NS_{0}=\sum_{n\ge 0}\C L_n +\sum_{r>0}\C G_r$.  Then every
non-trivial irreducible $\NS_+$-module is of the form ${\rm
Ind}_{\ns_{0}}^{\ns_+}\C_{\lambda}$, where $\C_\lambda$ is a one dimensional  
irreducible representation of $\NS_{0}$, on which $L_0$ acts as the scalar
$\lambda\in\C^*$ and $L_k$ and $G_r$ act trivially for $r,k>0$.

\proof We define a filtration on the Lie superalgebra $\NS_+$ as follows:
$\NS_i$ is the subalgebra
spanned by $\{L_{i\over 2},G_{{i+1}\over 2},L_{{i+2}\over 2},G_{{i+3}\over
2},\cdots\}$, if $i$ is even.  If $i$ is odd, then $\NS_i$ is spanned by the
linearly independent vectors $\{G_{i\over 2},L_{{i+1}\over 2},G_{{i+2}\over
2},L_{{i+3}\over 2},\cdots\}$.  We set $\dd=G_{-{1\over 2}}$ so that
$[\dd,\dd]=2L_{-1}$.  We then have $\NS_+=\C \dd+\C[\dd,\dd]+\NS_0$.  Hence  
by Theorem 3.1
every non-trivial irreducible conformal module over $\NS_+$ is of the form
${\rm Ind}_{\ns_{0}}^{\ns_+}U$, where $U$ is an irreducible $\NS_0$-module.
Now we use Lemma 3.2 and the result follows. $\qed$

Using similar filtrations as in Corollary 3.5 one proves Corollaries 3.6 and  
3.7 below.  Although never used, Corollary 3.6 is stated for the sake of  
completeness.

\proclaim{Corollary 3.6.}  Let $\G$ be a direct sum of finite-dimensional
simple Lie superalgebras. Let $\L=(\C{\partial\over {\partial
t}}+\C(\theta{\partial\over {\partial t}}-{\partial\over
{\partial\theta}}))\lsemi\G[t,\theta]$.  Then
every irreducible non-trivial conformal module of $\L$ is either of the form
${\rm Ind}_{\g[t,\theta]}^{\l} U$, where $U$ is a finite-dimensional irreducible
representation of $\G$ or it is an irreducible two dimensional
representation of the subalgebra $\C{\partial\over{\partial
t}}+\C({\theta{\partial\over{\partial t}}-{\partial\over{\partial
\theta}}})$, on which ${\partial\over{\partial t}}$ acts as a non-zero scalar
and $\G[t,\theta]$ acts trivially.

\proclaim{Corollary 3.7.} Let $\L=\NS_+\lsemi\G[t,\theta]$.  Set
$\L_0=\NS_0\lsemi\G[t,\theta]$.  Then every non-trivial
conformal module over $\L$ is of the form ${\rm Ind}_{\l_0}^{\l}U$, where
$U$ is a finite-dimensional non-trivial $\G\oplus\C L_0$-module, on which
$L_k$, $G_r$ for all $k,r>0$ and $\G[t]t+\G[t]\theta$ act trivially.

\noindent{\it Remark 3.2.} Corollary 3.7 gives a $3$-parameter family of
irreducible conformal modules over $\NS\lsemi\tilde{\G}_{\rm super}$.  As  
before we will denote the modules over the corresponding conformal  
superalgebra $R(\NS)\lsemi R(\tilde{\G}_{\rm super})$
by $M_{\ns\lsemi\tilde{\g}_{\rm super}}(\Lambda,\alpha,\Delta)$, where
$\Lambda$ is a dominant integral weight of $\G$, $\alpha,\Delta\in\C$.  The  
conditions for irreducibility of this module is as in Remark 3.2.

Translating the above corollaries into the language of modules over
conformal superalgebras (using Proposition 2.1 and Remark 2.2) we obtain
the following

\proclaim{Theorem 3.2.} Let $R(\V)$, $R(\tilde{\G})$, $R(\NS)$ and
$R(\tilde{\G}_{\rm
super})$  stand for the Virasoro, the current, the Neveu-Schwarz and the
super current conformal (super)algebras, respectively.  Let $R(\V)\lsemi
R(\tilde{\G})$ and
$R(\NS) \lsemi R(\tilde{\G}_{\rm super})$ denote their respective
semidirect sums.
Then the following is a complete list of their finite non-trivial
irreducible modules:
\item{1.} $M_{\v}(\alpha,\Delta)$ and $M_{\ns}(\alpha,\Delta)$, where
$\alpha,\Delta\in\C$ with $\Delta\not=0$.
\item{2.} $M_{\tilde{\g}}(\Lambda)$ and $M_{\tilde{\g}_{\rm
super}}(\Lambda)$, where $\Lambda$ is a non-zero dominant integral weight,  
or their quotients.
\item{3.} $M_{\v\lsemi\tilde{\g}}(\Lambda,\alpha,\Delta)$ and
$M_{\ns\lsemi\tilde{\g}_{\rm super}}(\Lambda,\alpha,\Delta)$, where $\Lambda$  
is a non-zero dominant integral weight and $\alpha,\Delta\in\C$ or else if  
$\Lambda=0$, then $\Delta\not=0$.

It was shown in [2] that every semisimple finite
conformal algebra is a direct sum of conformal algebras of the form $R(\V)$,
$R(\tilde{\G})$ and $R(\V)\lsemi R(\tilde{\G})$. The results of this section give
a  description of all finite irreducible modules over all finite
semisimple conformal algebras.  Namely we have the following

\proclaim{Proposition 3.1.} Let $R=R_1\oplus R_2\oplus \cdots \oplus R_n$
be a finite semisimple conformal algebra.  Suppose that
$R_i$ is either $R(\V)$ or $R(\V)\lsemi R(\tilde{\G})$ for some $i$.  Let $V$ be a
finite irreducible conformal module of $R$ whose restriction to $R_i$
is non-trivial.  Then the restriction of $V$ to all
$R_j$ is trivial for $i\not=j$.

\proof Since $R_i$ is either the conformal algebra $R(\tilde{\G})$ or  
$R(\V)\lsemi R(\tilde{\G})$, there exists $L_{(0)}^i\in R_i$ such that  
$[\dd-L_{(0)}^i,R_i]=0$.
Choose any $k\not= i$ and consider irreducible modules over the conformal  
algebra $R_i\oplus R_k$.  By Proposition 2.1 we are to consider modules over  
$\C\dd\lsemi((R_i)_+\oplus (R_k)_+)\cong
(R_i)_+\oplus(\C(\dd-L_{(0)}^i)\lsemi (R_k)_+)$.  Now a non-trivial
irreducible $(R_i)_+$-module is free over $\C[L_{(0)}^i]$ and a
non-trivial irreducible $\C(\dd-L_{(0)}^i)\lsemi (R_k)_+$-module is free
over $\C[\dd-L_{(0)}^i]$ by above discussion.  Thus their tensor product is  
free over $\C[L_{(0)}^i]\otimes\C[\dd-L_{(0)}^i]$ and hence cannot be finite  
over $\C[\dd]$. $\qed$

Proposition 3.1 and Theorem 3.2 imply immediately

\proclaim{Theorem 3.3.} Let $R=R_1\oplus R_2\oplus \cdots \oplus R_n$
be a finite semisimple conformal algebra, where each $R_i$ is either simple  
or of the form $R(\V)\lsemi R(\tilde{\G})$.  Then $R$ has a faithful finite
module if and only if either all $R_i$ are current conformal
algebras or $n=1$.

\bigskip
\bigskip
\centerline  {\bf REFERENCES }
\bigskip
\frenchspacing
\medskip

\item{1.} Cheng, S.-J.; Kac, V.~G.; Wakimoto, M.: Extensions of conformal
modules, preprint, 1997.

\item{2.} D'Andrea, A.; Kac, V.~G.: Structure theory of finite conformal algebras.

\item{3.} Kac, V.~G.: Vertex algebras for beginners. Providence: AMS
University lecture notes vol.~{\bf 10} 1996.

\item{4.} Kac, V.~G.: The idea of locality, in Proceedings of the XXI International colloquium on group theoretical methods in physics, Goslar, Germany, 1996.

\item{5.} Kac, V.~G.; Rudakov, A.~N.: Representations of simple finite conformal
superalgebras.
 
\item{6.} Rudakov, A.~N.: Irreducible representations of infinite-dimensional
Lie algebras of Cartan type, Math. USSR-Izvestija 8, 836-866 (1974).

\end